\title{Virtual algebraic Lie theory: 
  Tilting modules and Ringel duals for blob algebras}
\author{P P Martin and S Ryom--Hansen \\ \myaddress}
\documentclass[12pt]{article} 
\usepackage{epic} 
\usepackage{eepic} 
\usepackage{latexsym}
\usepackage{amssymb}
\usepackage[all]{xy}
\providecommand{\noglossaryignore}[1]{}
\newcommand{\globalglossaryentry}[3]{\makebox[1.5in][l]{\tt $\backslash${#1}} 
\makebox[1.1in][l]{{$#2$}} \makebox[2.5in][l]{{#3}}\newline} 
\newcommand{\newcommandabbreviation}[3]{\newcommand{#1}{#2}%
\noglossaryignore{\globalglossaryentry{#3}{#2}{}}}
\newcommand{\renewcommandabbreviation}[3]{\renewcommand{#1}{#2}%
\noglossaryignore{\globalglossaryentry{#3}{#2}{}}}
\newcommand{\newcommandmacro}[4]{\newcommand{#1}{#2}%
\noglossaryignore{\globalglossaryentry{#3}{#2}{#4}}}
\newcommand{\gge}[3]{\noglossaryignore{\globalglossaryentry{#1}{#2}{#3}}}
\newcommand{\myaddress}%
{\parbox{3in}{\footnotesize \begin{center} 
Mathematics Department, City University, \\  
Northampton Square, London EC1V 0HB, UK.\end{center}}}
    
\newcounter{minidef}[section]

\newcounter{minicapt}

\newtheorem{theo}{Theorem}      \newtheorem{cor}{Corollary}[theo]   
\newtheorem{de}{Definition}     \newtheorem{pr}{Proposition} 
\newtheorem{co}{Corollary}[pr]  \newtheorem{rem}{Remark} 
\newtheorem{lem}{Lemma} 
\noglossaryignore{GREEK ETC.\newline}
\newcommandabbreviation{\e}{\epsilon}{e}        
\newcommandabbreviation{\lam}{\lambda}{lam}  
\newcommandabbreviation{\la}{\langle}{la}        
\newcommandabbreviation{\ran}{\rangle}{ran}
\newcommandabbreviation{\ha}{\#}{ha}             
\newcommandabbreviation{\rmap}{\rightarrow}{rmap}
\newcommandabbreviation{\aaa}{\alpha}{aaa}        
\newcommandabbreviation{\ab}{\alpha,\beta}{ab}
\newcommandabbreviation{\aab}{a(\ab )}{aab}       
\noglossaryignore{\newline RINGS\newline}
\newcommandabbreviation{\HH}{H \!\!\! I}{HH}               
\newcommandabbreviation{\C}{\mathbb C}{C}
\newcommandabbreviation{\N}{\mathbb N}{N}   
\newcommandabbreviation{\Z}{\mathbb Z}{Z}      
\renewcommandabbreviation{\Re}{\mathbb R}{Re}
\newcommandabbreviation{\R}{{\mathbb R}}{R}
\newcommandabbreviation{\Q}{\mathbb Q }{Q}
\renewcommandabbreviation{\H}{\mathbb H }{H}
\noglossaryignore{\newline SYMMETRIC GROUP\newline}
\def\Sym(#1){\Sigma(#1)}                   
\gge{Sym(-)}{\Sym(-)}{symmetric group on - objects}
\def\Sy(#1){\Sigma_{#1}}                   
\gge{Sy(-)}{\Sy(-)}{symmetric group irreducible -}
\def\sym(#1){\mbox{\LARGE s}(#1)}        
\gge{sym(-)}{\sym(-)}{symmetric group on - objects (variant)}
\def\sy(#1){\mbox{\LARGE s}({#1})}        
\gge{sy(-)}{\sy(-)}{symmetric group irreducible - (variant)}
\newcommandmacro{\cs}{\C \, \sy(n)}{cs}{symmetric group algebra over $\C$}
\noglossaryignore{\newline PARTITIONS/SETS\newline}
\newcommand{\Nset}[1]{\underline{#1}}
\gge{Nset\{-\}}{\Nset{-}}{set of natural numbers to -} 
\def\nset(#1){ \{ #1 \}_{ \underline{n} }} 
\gge{nset(-)}{\nset(-)}{a set $-\times\Nset$}
\def\ul(#1){_{\underline{#1}}}             
\gge{ul(-)}{{}\ul(-)}{subscript underline -}
\def\Ee(#1){{\bf E}_{#1}}                  
\gge{Ee(-)}{\Ee(-)}{set of equivalence relations on set -}
\def\Eee(#1){{\bf E}_{\{ #1 \}_{\underline{n}}}}   
\gge{Eee(-)}{\Eee(-)}{ditto for nset}
\def\Een(#1,#2){{\bf E}_{\{ #1 \}_{\underline{#2}}}}   
\def\Ssn(#1,#2){{\bf S}_{\{ #1 \}_{\underline{#2}}}}   
\def\Ss(#1){{\bf S}_{#1}}                  
\def\Sss(#1){{\bf S}_{\{ #1 \}_{\underline{n}}}}   
\def\bbc(#1){((\beta_1)(\beta_2)...(\beta_{#1}))}      
\newcommandmacro{\Ln}{{\Gamma}^{n}}{Ln}{large index set}
\newcommandmacro{\LnQ}{{\Gamma}^{n}_Q}{LnQ}{index set}
\newcommandmacro{\Zz}{\zeta}{Zz}{`shape' function}
\noglossaryignore{\newline PARTITION ALGEBRA\newline}
\def\ka(#1){\kappa_{#1}}                   
\def\Sm(#1){\Sigma_{#1}}                   
\newcommandmacro{\com}{\bullet}{com}{bullet composition}
\newcommandmacro{\enm}{\; e^n(\! m\! ) \;}{enm}{product of idempotents}
\def\Ai(#1){ A^{ #1 \cdot } }              
\def\Aij(#1,#2){ A^{ #1  #2 } }            
\newcommandmacro{\One}{\mbox{\bf $1 \!\!\! 1$}}{One}{algebra unit 1}
\newcommandmacro{\Bp}{B_p}{Bp}{partition basis}
\def\Bb(#1){B_p[#1]}                       
\def\Pp(#1){P_n[#1]}                       
\def\Ps(#1){P_n[#1] \! /}                  
\newcommandmacro{\Ph}{\hat{P}}{Ph}{P hat  algebra}
\def\Is(#1){\sim^{#1}}                     
\noglossaryignore{\newline MODULES\newline}
\def\Wm(#1){{\cal S}_{#1}}                 
\gge{Wm(-)}{\Wm(-)}{Weyl module with index -}
\def\wm(#1,#2){{}_{#1}{\cal S}_{#2}}       
\gge{wm(-1,-)}{\wm(-1,-)}{Weyl module with index -}
\def\Ind(#1,#2,#3){\mbox{Ind}_{#1}^{#2}#3} 
\gge{Ind(-1,-2,-)}{\Ind(-1,-2,-)}{induction}
\def\Res(#1,#2,#3){\mbox{Res}_{#1}^{#2}#3} 
\gge{Res(-1,-2,-)}{\Res(-1,-2,-)}{restriction}
\newcommandabbreviation{\weyl}{standard}{weyl}
\newcommandabbreviation{\mod}{\mbox{mod}}{mod}
\newcommandabbreviation{\head}{\mbox{head }}{head}
\newcommandabbreviation{\Weyl}{Weyl}{Weyl}
\def\SS(#1){{\cal S}_{#1}}                 
\gge{SS(-)}{\SS(-)}{Specht/Weyl module index -}
\def\LL(#1){{\cal L}_{#1}}                 
\gge{LL(-)}{\LL(-)}{simple module index -}
\noglossaryignore{\newline FUNCTORS/MAPS\newline}
\newcommandmacro{\Gg}{{\cal G}}{Gg}{G Functor}
\newcommandmacro{\Fg}{{\cal F}}{Fg}{F Functor}
\newcommandmacro{\ra}{\rightarrow}{ra}{}
\def\ses(#1,#2,#3){0\ra #1 \ra #2 \ra #3 \ra 0}   
\gge{ses(1,2,-)}{\ses(1,2,-)}{\hspace{.5in} short exact sequence}
\def\starr(#1){ \stackrel{ #1 }{\longrightarrow} }
\gge{starr(-)}{\starr(-)}{}
\newcommandmacro{\doublerightarrow}{\; -\!\!\! -\!\!\!\!\!\! \gg \;}
{doublerightarrow}{}
\noglossaryignore{\newline PARTITION ALGEBRA MAPS\newline}
\newcommandmacro{\smap}{s}{smap}{`inclusion' map}
\newcommandmacro{\tmap}{t}{tmap}{$ P_n -> S_n$}
\newcommandmacro{\pmap}{\psi}{pmap}{$ S_n -> P_n $}
\noglossaryignore{\newline MISC.\newline}
\def\Amap(#1){{\cal A}_{#1}}               
\gge{Amap(-)}{\Amap(-)}{}
\def\Rr(#1){R_{#1}}                        
\gge{Rr(-)}{\Rr(-)}{restriction of E}
\def\Cr(#1){C_{#1}}                        
\gge{Cr(-)}{\Cr(-)}{restriction of E to N}
\newcommandmacro{\Tm}{{\cal T}}{Tm}{Transfer Matrix}
\def\On(#1){{\cal I}_{#1}}
\gge{On(-)}{\On(-)}{}
\newcommandmacro{\UU}{\underline{\sqcup}}{UU}{}  
\newcommandmacro{\UUU}{\sqcup}{UUU}{}  
\newcommandmacro{\Vq}{V_Q^{\otimes n}}{Vq}{Potts config. space}
\def\bs(#1,#2){\mbox{{\Large $\ast$}}^{#1}_{#2}}  
\gge{bs(-,-)}{\bs(-,-)}{general plumbing multiplier}
\newcommand{\ignore}[1]{}
\gge{ignore\{-\}}{\ignore{-}}{ignore argument!}
\def\choo(#1,#2){ \left( \begin{array}{c} #1 \\ #2 \end{array} \right) } 
\gge{choo(-1,-)}{\choo(-1,-)}{choose}
\newcommand{\Qed}{$\Box$}
\gge{Qed}{\mbox{\Qed}}{QED}
\def\staq(#1){\stackrel{#1}{=}}            
\gge{staq(-)}{\staq(-)}{}
\def\stam(#1){\stackrel{#1}{\rightarrow}}  
\gge{stam(-)}{\stam(-)}{}
\def\mat{ \left( \begin{array} }    
\def\tam{ \end{array}  \right) }
\gge{mat/tam}{...}{matrix delimiters}
\newcommand{\beq}{\begin{equation} }
\def\eql(#1){ \begin{equation} \label{#1} 
}
\newcommand{\eq}{\end{equation} }
\def\eqal(#1){\begin{eqnarray} \label{#1} }
\def\eqa{\end{eqnarray} }
\def\lab(#1){\label{#1}
}
\def\prl(#1){ \begin{pr} \label{#1} 
}
\def\del(#1){ \begin{de} \label{#1} 
}
\gge{smeq\{-\}}{...}{small equation}
\gge{fneq\{-\}}{...}{very small equation}
\noglossaryignore{\newline HECKE/BLOB\newline}
\newcommandmacro{\Hnq}{H_n(q)}{Hnq}{ * freestanding symbol}
\newcommandmacro{\Hn}{H_n}{Hn}{      *-mod etc.}
\newcommandmacro{\A}{{\cal A}}{A}{}
\newcommandmacro{\Cwts}{C}{Cwts}{}
\newcommandmacro{\CA}{{\cal A}}{CA}{}

\newcommandmacro{\calA}{{\cal A}}{calA}{}
\newcommandmacro{\modi}{\mbox{Mod} }{modi}{was mod not modi!}
\newcommandmacro{\Wgen}{{\Bbb S}}{Wgen}{}
\def\ol(#1){\overline{#1}}
\newcommandmacro{\st}{\mbox{St}}{st}{}
\def\CMult(#1,#2){(#1:#2)}
\def\CM(#1,#2){( #1 : #2 )}
\def\FMult#1,#2{(#1:#2)}
\def\CF#1,#2{(#1:#2)}

\newcommandmacro{\Top}{\mbox{Top}}{Top}{}
\newcommandmacro{\Soc}{\mbox{Soc}}{Soc}{}
\newcommandmacro{\Head}{\mbox{Head}}{Head}{}
\newcommandmacro{\Filt}{{\cal F}}{Filt}{}
\newcommandmacro{\Mod}{\mbox{mod}}{Mod}{}
\newcommandmacro{\Resi}{\mbox{Res }}{Resi}{was without i!}
\newcommandmacro{\Indi}{\mbox{Ind }}{Indi}{was without i!}
\def\RR(#1,#2){R^{#1}_{#2}}   
\def\TT(#1,#2){T^{#1}_{#2}}   
\def\implies{\Rightarrow}

\newcommandmacro{\Ann}{\mbox{Ann}}{Ann}{}
\newcommandmacro{\Cen}{\mbox{Cen}}{Cen}{}
\newcommandmacro{\End}{\mbox{End}}{End}{}
\newcommandabbreviation{\semisimple}{semisimple}{semisimple}
\newcommandabbreviation{\Bratteli}{Bratteli}{Bratteli}
\newcommandabbreviation{\JBC}{Jones Basic Construction}{JBC}
\newcommandabbreviation{\pa}{partition algebra}{pa}
\newcommandabbreviation{\TM}{transfer matrix}{TM}
\newcommandabbreviation{\PM}{Potts model}{PM}
\newcommandabbreviation{\QSC}{quantum spin chain}{QSC}
\newcommandabbreviation{\Hamiltonian}{Hamiltonian}{Hamiltonian}
\newcommandabbreviation{\YS}{Young symmetrizer}{YS}

\newcommand{\beqa}{\begin{eqnarray}}%
\newcommand{\eeqa}{\end{eqnarray}}%

\newcommand{\vardelta}{\delta_e}%

\newcommand{\mU}{{\cal  U}}
\newcommand{\ym}{{m}}
\newcommand{\TL}{ Temperley--Lieb}
\newcommand{\qh}{quasihereditary}

\newcommand{\chix}{\chi}
\newcommand{\egen}{e}
\newcommand{\Uqsl}{U_q{\frak s \frak l}}
\newcommand{\muM}{\mu}
\newcommand{\qr}{q}%
\newcommand{\lex}{lexicographic} 

\newcommand{\seq}{\mbox{seq}}%
\newcommand{\res}{\mbox{Res}}%

\newcommand{\SSv}{{\cal S}_n}%
\newcommand{\hash}{\#}%

\newcommand{\vargamma}{\gamma}%
\newcommand{\RK}{R}%
\newcommand{\Permmod}{M}%
\newcommand{\permmod}[2]{\Permmod_{#2}(#1)}%
\newcommand{\standard}[2]{\Delta_{#2}(#1)}%
\newcommand{\tilting}[2]{T_{#2}(#1)}%
\newcommand{\ourrep}{\rho}
\newcommand{\altrep}{\ourrep'}
\newcommand{\our}{}%
\newcommand{\altM}{{\cal M}}%
\newcommand{\SCh}[1]{(#1)}%

\newcommand{\idem}{\epsilon}
\newcommand{\signU}{+}
\newcommand{\hgt}{\triangleright}
\newcommand{\hlt}{\triangleleft}

\begin{document} \maketitle
 \newcommand{\ignoreifnotdraft}[1]{
}
\ignoreifnotdraft{
\pagestyle{myheadings} \markboth{Draft}{\today}
}

\ignoreifnotdraft{
\noindent
\begin{tabular}{l} 
  {\tiny \filename (Draft)} \hspace{4.4in} Jan 1997 \\ \hline 
\end{tabular}
}
\ignoreifnotdraft{ \newpage } 

\section{Introduction}%

In this paper we construct a representation of the blob algebra
\cite{MartinSaleur94a} 
over a ring allowing base change to every interesting
(i.e. non--semisimple) specialisation which, in \qh\ 
specialisations, passes to a full tilting module. 

The Temperley--Lieb algebras are  a tower 
$T_0(q) \subset T_1(q) \subset ..$
of one--parameter finite dimensional algebras \cite{TemperleyLieb71}, 
each with a basis independent of $q$. 
These algebras are \qh\ 
\cite{ClineParshallScott88,DlabRingel89} 
except in case $q+q^{-1}=0$. 
Accordingly one may in principle construct tilting modules, full
tilting modules, and corresponding Ringel duals. 
In fact, if $V$ is a free module of
rank 2 over the ground ring then $T_n(q)$ has an action on 
$V^{\otimes  n}$, and it is straightforward to show (see later) that 
$V^{\otimes  n}$ is a full tilting module 
in the \qh\ cases.
Since $V^{\otimes  n} $ exists over the ground ring, 
the Ringel dual can be constructed without having
to pick a specialisation.
The cases of $n$ finite of this dual
are a nested sequence of quotients of the quantum group $\Uqsl_2$ 
\cite{Jimbo85,Drinfeld85}. 
This $q$--deformable duality and glorious limit structure \cite{Jantzen}
(more usually observed with $\Uqsl_2$ as the starting point) 
provides the mechanism for massive exchange of representation theoretic
information between the two sides 
\cite{Ringel91,Donkin98,ErdmannHenke02,CoxGrahamMartin01,KleshchevSheth99}. 
In particular the weight theory of $\Uqsl_2$ controls the
representation theory of $T_n(q)$ for all $n$ simultaneously (as
localisations of a global limit). 


The blob algebras are  a tower $b_0 \subset b_1 \subset ..$ of
two--parameter finite dimensional algebras (and $b_n \supset
T_n(q)$). They are \qh\ except 
at a finite set of parameter values.
Accordingly one may in principle construct tilting modules and so on. 
{\em Ab initio} one would have to expect such a construction to depend
on the specialisation, 
as {\em indecomposable } tilting modules do \cite{MartinWoodcock2000}. 
On the other hand, it turns out \cite{MartinWoodcock02} that 
$b_n$ has an action on 
$V^{\otimes  2n}$, and in this paper we show that
$V^{\otimes  2n}$ is a full tilting module 
in the \qh\ cases.


Historically, $T_n(q)$ and $\Uqsl_2 $  
were studied extensively separately,
before the full tilting module/Ringel duality connection was known, 
but if one side, and the appropriate full tilting module, had been discovered
first, the passage to the Ringel dual would rightly have been
regarded as quite a significant spin--off!
The $b_n$ tilting property of $V^{\otimes  2n}$ is a striking result, 
in as much as it places us in a position analogous to
this
(as it were, before the discovery of quantum groups). 


Suitably prepared, the blob algebra may be regarded as a quotient of
the Ariki--Koike algebra, which itself is a quotient of the 
affine Hecke algebra \cite{ArikiKoike94,Hoefsmit74}. 
Thus the representation theory of the blob algebra is 
part of the representation theory of the Ariki--Koike and of the affine 
Hecke algebra, this last point being the basic idea of 
\cite{MartinSaleur93,MartinSaleur94a} (see also Graham and Lehrer's 
analysis \cite{GrahamLehrer01}).
Although $V^{\otimes  2n}$ can be regarded as a module of 
the Ariki--Koike algebra, it is {\em not} in any obvious way a sum of
permutation modules in the usual Ariki--Koike sense 
\cite{DipperJamesMurphy95,DipperJamesMathas97b,DipperJamesMathas97c}.  
In the absense of a {\em natural} tensor space 
(cf. \cite{ArikiTerasomaYamada95,SakamotoShoji99}), 
these permutation modules 
form the starting point for most  tilting related approaches to Ariki--Koike
representation theory.
Our approach is of an essentially different nature. 
In particular it gives a weight theory (in the sense mentioned above)
for $b_n$ which, for Ariki--Koike,  
would imply a structure
which it seems very unlikely to possess 
(see \cite{MartinWoodcock02}). 


The blob algebra, and certain generalisations, have been observed to
manifest several indicators of an underlying structure evocative
of algebraic Lie theory (such as the role played in their representation
theory by alcove geometry --- see \cite{MartinWoodcock02}). 
We wish to understand the underlying reasons for the extra structure. 
The \TL\ paradigm suggests that an appropriately prepared Ringel
dual is a good place to look (hence `{\em virtual} algebraic Lie theory'). 
For example, the results of 
\cite{MartinWoodcock2000} suggest that this ``dual blob algebra'' should be 
reminiscent of
the Kac--Moody quantum algebra $ U_q \hat{ \frak s \frak l}_2$.


Given the full background, 
a  natural approach to proving that $V^{\otimes  n}$ is a tilting
module for $T_n(q)$ is to use the duality itself 
(that this module is
a tilting module on the dual side is a direct consequence of the
general machinery of Donkin \cite{Donkin93,Donkin98} et al \cite{Erdmann94}). 
The challenge here is that for the blob {\em no such general machinery
  yet exists}. 
Accordingly
we include here a  proof in the \TL\ case which does not appeal to the 
algebraic Lie theory 
machinery, but only to quasiheredity. 
This exercise is motivated only by the need to understand how
such a proof might work, for use in the blob case. 
The blob case is then the main object of this paper 
(see section~\ref{blobtilt}). 


\subsection{Preliminaries}\label{Preliminaries}


The blob algebra $b_n$ is usually defined in terms of a certain basis
of diagrams and their compositions 
\cite{MartinSaleur94a},  
from which it derives its name. 
This blob algebra is isomorphic to an algebra defined by a
presentation \cite{CoxGrahamMartin01}.  
We will only need the presentation. 


For $\RK$ a commutative ring, 
$x$ an invertible element in $\RK$,
$q=x^2$, and $\vargamma,\vardelta \in \RK$, define
$b_n^{\RK}$
to be the $\RK$--algebra with 
generators $\{ 1,\egen,U_1,...,U_{n-1} \}$ 
and relations
\begin{eqnarray} 
 U_i U_i &=& (q+q^{-1}) U_i \label{TL001} \\
 U_i U_{i\pm 1} U_i &=& U_i \\
 U_i U_j &=& U_j U_i \hspace{1in} \mbox{($|i-j|\neq 1$)}  \label{TL003}
\end{eqnarray}
\begin{eqnarray*} 
U_1 \egen U_1 &=&  \vargamma U_1 \\
  \egen \egen &=& \vardelta \egen \\
    U_i \egen &=& \egen U_i   \hspace{1in} \mbox{($i \neq 1$)} . 
\end{eqnarray*}
By the isomorphism with the diagram algebra 
this algebra is a free $\RK$--module 
(with basis most conveniently described in terms of diagrams, 
however we do not otherwise need this basis, so we will not recall it
here --- see \cite{MartinSaleur94a}). 


It will be evident that $\egen$ can be rescaled to change $\vargamma$ and
$\vardelta$ by the same factor. 
Thus, if we require that $\vardelta$ is invertible, then we might as well
replace it by 1. 
(This brings us to the original two--parameter definition of the algebra.) 


For $k$ a field which is a $\RK$--algebra define 
$b_n = k \otimes_{\RK} b_n^{\RK}$. 
It is known \cite{MartinSaleur94a} that the representation theory of
$b_n$ falls into one of three distinct categories, 
depending on the number of integer values of $a$ for which 
$$ \vargamma [a]_q = \vardelta [a-1]_q $$
(where $[a]_q$ is the usual $q$--number). 
If there  is no solution then $b_n$ is semisimple 
and has trivial tilting theory. 
Accordingly it is convenient to reparameterize into
the following form: 
\eql(our_pres)  
\vargamma=q^{\ym-1} - q^{-\ym+1} , 
\hspace{.8cm}   
\vardelta=q^{\ym} - q^{-\ym} . 
\eq
Note  
that provided $m$ is integer (which includes all the
interesting cases) this parameterization has a lattice in $b_n^{\Z[q,q^{-1}]}$ 
(i.e. $\vargamma,\vardelta $  lie in $\Z[q,q^{-1}] $). 


Previously used parameterizations include
\newline
$\vargamma=\frac{[\ym-1]}{[\ym]}$, 
$\;$ $\vardelta=1$;
\newline
and
\newline
$\vargamma= \pm [\ym -1]$, 
$\;$ $\vardelta= \pm [\ym]$
\newline (see \cite{MartinSaleur94a}, \cite{MartinWoodcock2000}, 
\cite{CoxGrahamMartin01} respectively). 
Our form has the mild disadvantage that it is not a simple rescaling in
case $q-q^{-1} = 0$.

\section{The `crypto--tensor' representations}\label{the rep}

We now recall
the representations 
$\rho_{\chix}$ 
of $b_n$ defined in \cite[\S6.1]{MartinWoodcock02}. 

Set
$$\mU^q(\chix)=\mat{cccc} 
0&0&0&0\\ 
0&q&1&0 \\
0&1&q^{-1}&0 \\
0&0&0&\chix \tam$$
and $\mU^q =\mU^q(0)$. 


Let $V= \mbox{\rm span} \{v_1,v_2 \} $. 
Let $\seq_{}\{1,2\}$ denote the set of words of finite length in $\{1,2\}$, 
$\seq_n\{1,2\}$ the subset of words of length $n$,
and $\seq_n^{r}\{1,2\}$ the subset of this in which the number of 1s is $r$. 
For $w \in \seq_n\{1,2\}$ let $\hash^1(w)$ denote the number of 1s in
$w$ (so $w \in \seq_n^{r}\{1,2\}$ implies $\hash^1(w)=r$).  
Then $V^{\otimes n}$ has basis 
$\{ v_{i_1} \otimes v_{i_2} \otimes .. \otimes v_{i_n} \; | \; 
  i_1 i_2 .. i_n \in \seq_n\{1,2\} \}$. 
We will adopt the shorthand of writing the sequence for the basis
element. 
We ascribe the usual 
\lex\ 
order to this basis (11,12,21,22 and so on). 
Let $\mU^{q}(\chix)$ act on $V \otimes V$ 
with respect to this ordering of the basis. 

Let $\muM^{\qr,\chix}(U_i) \in \End(V^{\otimes n})$ 
be a matrix acting trivially on every tensor
factor except the $i^{th}$ and $(i+1)^{th}$, where it acts as
$\signU\mU^{\qr}(\chix)$. 
Write $\muM^{\qr}(U_i)$ for  $\muM^{\qr,0}(U_i)$. 
The \TL\ algebra $T_n(q)$ is the subalgebra of $b_n$ with generators 
$\{ 1,U_1,...,U_{n-1} \}$. 
The {\em tensor space} representation of $T_n(q)$ is given by $\muM^{q}$. 


Note that 
\eql(rst) 
(\mU^{s}\otimes \mU^{t}) 
   (1 \otimes \mU^{r}(\chix) \otimes 1) 
      (\mU^{s}\otimes\mU^{t})
= 
\left( \frac{r}{st} +\frac{st}{r} +\chix \frac{t}{s} \right)  
  (\mU^{s}\otimes \mU^{t})
\eq
for any $r,s,t,\chix$ (an explicit calculation). 




Suppose 
henceforth that 
there is an element $a \in K$ such that $a^4=-1$. 
Then $a^2 + a^{-2} = 0$. 
Set 
$$r= a^2 q^{\ym}$$ 
$$s  
 = a^5 x $$ 
$$t  
 = a^3 x $$ 
We have 
$$r+r^{-1} = a^2 ( q^{\ym} - q^{-\ym} ) $$
$$ s + s^{-1} = a^5 x + a^3 x^{-1} $$
$$ t + t^{-1} = a^3 x + a^5 x^{-1} $$
$$st=q$$
$$[2]_s [2]_t = [2]_q$$
$$ \frac{st}{r}+\frac{r}{st} 
   = a^2(q^{\ym -1} - q^{1-\ym}) $$ 

Then by equation~(\ref{rst}) there is an algebra homomorphism 
$$
\ourrep : b_n^{\Z[q,q^{-1}]}(q,\ym) 
  \longrightarrow End_{\Z[a,x,x^{-1}]}(V^{\otimes 2n})$$ 
given by 
\beqa \label{mape}
\ourrep : \egen & \mapsto & 
                        a^{-2} \muM^r(U_n)
\\
\ourrep: U_i  & \mapsto &  \muM^s(U_{n-i}) \muM^t(U_{n+i})
\eeqa
for $b_n$ in the  
form described in equation~(\ref{our_pres}). 


There is another algebra homomorphism  
$\altrep$ 
defined in exactly the same way
except that 
\beqa \label{mapex} 
\altrep: \egen & \mapsto & 
                        a^{-2} \muM^{r,[2]_r}(U_n)
\eeqa


Note that the $T_n(q)$--module 
$\muM^q$ has   {\em manifest} direct summands with basis
$\seq_n^{r}\{1,2\}$ ($r=0,1,..,n$), called {\em permutation modules}. 
Similarly  
$\ourrep$ and $\altrep$ 
have {\em manifest} 
direct summands with basis $\seq_{2n}^{r}\{1,2\}$ which we will again
call permutation modules. 
\footnote{
It should be emphasised that 
although $b_n$ is a quotient of the Hecke algebra of type--B (which 
itself is an Ariki--Koike algebra), the above permutation modules 
do {\bf not} coincide with the type--B or Ariki--Koike permutation modules
described in \cite{DipperJamesMurphy95,Mathas02x}: 
in the notation of \cite{MartinWoodcock02} 
the quotient map sends $ g_i+q $ to $U_i$, and so
$U_i$ will map a typical Ariki--Koike permutation module basis vector 
to a linear combination of precisely two basis vectors 
(see \cite{DipperJamesMurphy95}), 
which is clearly not the case in our situation. 
As a matter of fact, an Ariki--Koike permutation module is 
typically not a module for the blob algebra, even if its leading Specht
factor is one. 
}


Similarly evidently we have
\prl(Vres) The following are {\em manifest} direct sums
(i.e. respecting the basis):
\beqa \label{TLVres}
\res^{T_n}_{T_{n-1}} \muM^q &=& \muM^q \oplus \muM^q
\eeqa
\beqa \label{VVres}
\res^{b_n}_{b_{n-1}} \ourrep 
&=& \ourrep \oplus \ourrep \oplus \ourrep \oplus \ourrep
\eeqa
\end{pr}

\section{A \TL\ tilting module}



The goal of this section is to prove that the tensor space module
$V^{\otimes n } $ 
is a tilting module for the Temperley--Lieb algebra. This follows from
general results (for instance in \cite{Donkin98}), but we give
here a self--contained argument which later generalizes to the
blob algebra representation $\ourrep$.

\medskip


{ \em We will from now on assume that $ [2]_q \not= 0 $ over $k$}.
Then the Temperley-Lieb algebras $ T_n  = T_n(q)$ are \qh. 
In fact \cite{Martin91,MartinWoodcock2000}, 
setting $ \idem=\frac{1}{[2]_q} U_{n-1} $, we 
have that 
$ \idem $ is part of a heredity chain for $ T_n $, 
and since $\idem$ and $T_{n-2} \subset T_n$ commute, 
$T_n \idem$ is a right $T_{n-2}$--module, and indeed 
\eql(TL)
\idem T_n \idem = \idem T_{n-2}  \cong  T_{n-2} 
\eq

Let $ F $ be the localization functor 
$$ F: T_n \, \mbox{-mod} \rightarrow T_{n-2} \, 
\mbox{-mod}:\,\, M \mapsto \idem M $$
and let $ G $ be the globalization functor 
$$ G: T_{n-2} \,\mbox{-mod} \rightarrow T_{n} \,\mbox{-mod}:\,\, 
N \mapsto T_{n} \idem
\otimes_{ T_{n-2} } N 
$$
Note that $ F $ is exact and $ G $ is right exact, being the 
left adjoint of $ F $.

\medskip

\subsection{Homological considerations}

Since the categories $ T_n \, \mbox{--mod} $ are \qh\  they come
with standard modules $ \Delta_n(\lambda) $, costandards $ \nabla_n(\lambda) $,
simples $ L_n(\lambda) $, their projective covers $ P_n(\lambda ) $, 
injective envelopes $ I_n(\lambda ) $ and tiltings $ T_n(\lambda) $ 
for $ \lambda \in \Lambda_n =  \{n, n-2,   \ldots, 0/1  \} $.
(The heredity order $\hgt $ is the reverse of the natural order on
$\Lambda_n$ as a subset of $\Z$.) 
The following statements can be found in 
appendix~A1 of Donkin's book \cite{Donkin98} (in a much more 
general setting than ours):
\prl(Apply F)
Assume that $ \lambda \in \Lambda_{n-2} $. Then 

\noindent
i) $F L_n(\lambda) = L_{n-2}(\lambda) $. 

\noindent
ii) $F \Delta_n(\lambda) = \Delta_{n-2}(\lambda) $ and 
$F \,  \nabla_n(\lambda) = \nabla_{n-2}(\lambda) $.

\noindent
iii) $ F P_n(\lambda ) = P_{n-2}(\lambda )$ and $ F I_n(\lambda ) 
= I_{n-2}(\lambda ) $.

\noindent
Otherwise (i.e. for $ \lambda = n $) we have that $F L_n(\lambda) = 
F \Delta_n(\lambda) = F \nabla_n(\lambda ) = 0$.
\end{pr}


Our next step is to investigate the application of $ G $ to these modules.
Write $ M \in {\cal F }_{n}(\Delta) $ if $M \in T_n\mbox{--mod}$ has a
standard filtration.  
Write $(M : \Delta_{n}(\mu))$ for the multiplicity of
$\Delta_{n}(\mu)$ as a filtration factor of $M$. 
We need the following Proposition:

\prl(apply G) 
If 
$ M \in {\cal F }_{n-2}(\Delta)\, $, 
then $GM$ also has a standard filtration. 
Furthermore the standard multiplicity is 
$$
( G M : \Delta_{n}(\mu) ) = \left\{
\begin{array}{ll}
( M : \Delta_{n-2}(\mu) ) &  \mbox{if  }  \mu \in \Lambda_{n-2} \\
0 & \mbox{otherwise}  \\
\end{array}
\right.
$$
\end{pr}
{\em Proof:}
By Donkin's homological criterion, see e.g. \cite[A2.2 (iii)]{Donkin98}, 
we know that $ G M \in {\cal F }_n(\Delta) $ 
if and only if 
$$ \mbox{ Ext}^1_{T_{n}} ( G M, \nabla_n(\mu ) )=0 
\,\,\,\,\, \forall  \mu \in \Lambda_{n} $$
Applying $ \mbox{ Hom}_{T_{n}} ( G M , - )  $ to the short exact 
sequence
$ \nabla_n(\mu ) \hookrightarrow I_n(\mu) \twoheadrightarrow Q_n(\mu) $
(defining $ Q_n(\mu)$) 
yields a long exact sequence whose first terms are
$$ \begin{array}{lll}
0  \rightarrow \mbox{ Hom}_{T_{n}} ( G M,\nabla_n(\mu ) )  \rightarrow &
\mbox{ Hom}_{T_{n}} ( G M,I_n(\mu ) )  \rightarrow & 
\mbox{ Hom}_{T_{n}} ( G M ,Q_n(\mu ) )  \\ \,\,\,\,\, \rightarrow
\mbox{ Ext}_{T_{n}}^1 ( G M ,\nabla_n(\mu ) )  \rightarrow & 0 &
\end{array} $$
Assume first that $ \mu \in \Lambda_{n-2} $. 
Then by  Proposition~\ref{Apply F} 
and adjointness, the first three terms of this sequence become
\eql(exact)
 0 \rightarrow \mbox{ Hom}_{T_{n-2}} ( M,\nabla_{n-2}(\mu ) ) \rightarrow
\mbox{ Hom}_{T_{n-2}} ( M,I_{n-2}(\mu ) ) \rightarrow 
\mbox{ Hom}_{T_{n-2}} ( M,F Q_n(\mu ) ) 
\eq
which coincides with the beginning of the long exact sequence that arises
from applying  
$ \mbox{ Hom}_{T_{n-2}} ( M, - )  $ to 
$ F \nabla_n(\mu ) \hookrightarrow F I_n(\mu) \twoheadrightarrow 
F Q_n(\mu) $. 
But then by (the easy direction of) the criterion, 
the last map of (\ref{exact}) is surjective and thus
$ \mbox{ Ext}_{T_{n}}^1 ( G M ,\nabla_n(\mu ) ) = 0 $ as claimed. 


Assume then that $ \mu = n $. When we apply $ F $ to the short exact 
sequence 
$  \nabla_n(n ) \hookrightarrow  I_n(n) \twoheadrightarrow 
 Q_n(n) $ we get an isomorphism
$  F I_n(n) \cong F Q_n(n) $. 
So applying 
$ \mbox{Hom}_{T_{n}}(G M, - ) $ to it, we get
$$ 
0 \rightarrow \mbox{ Hom}_{T_{n}} ( G M,\nabla_{n}(n ) ) \rightarrow
\mbox{ Hom}_{T_{n}} (G M,I_{n}(n ) ) \rightarrow 
\mbox{ Hom}_{T_{n}} ( G M, Q_n(n ) ) $$
where by adjointness the last map is an isomorphism, so 
also in this case 
$ \mbox{ Ext}_{T_{n}}^1 ( G M ,\nabla_n(\mu ) ) = 0 $, and 
the criterion applies.


To get the multiplicity statement, recall that since 
$ M  \in {\cal F}_{n-2}( \Delta) $, $ G M  \in {\cal F}_{n}( \Delta) $,
we have 
$$ 
( G  M : \Delta_n(\mu) ) =\mbox{dim Hom}_{T_{n}}(G M , \nabla_n(\mu) )= 
\mbox{dim Hom}_{T_{n-2}}(M,F \nabla_n(\mu) ) 
$$ 
which is zero for $ \mu = n $, while for $ \mu \in \Lambda_{n-2} $  
$$  
\mbox{dim Hom}_{T_{n-2}}(M , \nabla_{n-2}(\mu) )= 
( M : \Delta_{n-2}(\mu)) 
$$  
and we are done. \Qed
\begin{rem}
$ G $ does {\bf not} map $ { \cal F }_{n-2}(\nabla)\,\, $to 
$ { \cal F }_{n}(\nabla) $.
\end{rem}

We may now, incidentally, prove:
\begin{co}
\eql(dtool) 
L^{i}G \Delta_{n-2}(\lambda)  \,\, = \,\,
\left\{
\begin{array}{ll}
\Delta_n(\lambda) & \mbox{    if } i = 0 \\
0      & \mbox{    otherwise }
\end{array}
\right.
\eq 
\end{co}
 
{\em Proof:}
The $ i=0 $ case follows from Proposition~\ref{apply G}. By 
applying $ G $ to the short exact sequence 
$ K_{n-2}(\lambda) 
\hookrightarrow P_{n-2}(\lambda ) \twoheadrightarrow 
\Delta_{n-2}(\lambda )$, 
the last terms of the resulting long exact cohomology sequence are as follows:
$$ L^1 G \Delta_{n-2}(\lambda ) \hookrightarrow 
G K_{n-2}(\lambda) \rightarrow G P_{n-2}(\lambda ) \twoheadrightarrow G 
\Delta_{n-2}(\lambda )
$$
Now $ K_{n-2}(\lambda) $ has a standard filtration, so once again 
using Proposition~\ref{apply G} we get
$ [G P_{n-2}(\lambda )] = [G \Delta_{n-2}(\lambda )]+ 
[G K_{n-2}(\lambda) ] $, where as usual $[M]$ denotes the image of 
$M \in T_n \mbox{-mod}$ in the Grothendieck group. 
But then $ [L^1 G \Delta_{n-2}(\lambda)] = 0 $
and thus $ L^1 G \Delta_{n-2}(\lambda )= 0 $.

To get the vanishing of the higher
$ L^i G \Delta_{n-2}(\lambda ) $, we use induction from above 
(with respect to the heredity order) on $ \lambda $. 
If $ \lambda $ is maximal ($\lambda =0$ or 1) 
then $ P_n(\lambda) = \Delta_n(\lambda) $ 
and there is nothing to prove. Otherwise note that by the 
long exact sequence, 
$ L^{i} G K_{n-2}(\lambda ) = L^{i+1} G \Delta_{n-2}(\lambda ) $ for 
$ i > 0 $ so it is enough to show that $ L^{i} G K_{n-2}(\lambda ) = 0 $
for these $i$. 
But only $\Delta( \mu) $ with $ \mu \hgt \lambda $ occur 
in the standard filtration of $ K_{n-2} $
so the induction hypothesis applies to them. 
Let $\nu$ be such that $K_{n-2} \twoheadrightarrow \Delta(\nu)$. 
Considering the short  exact sequence 
$ K^{\nu}_{n-2} \hookrightarrow K_{n-2} 
\twoheadrightarrow \Delta(\nu) $ (defining $ K^{\nu}_{n-2}$), 
we get that
$ L^i G K^{\prime}_{n-2} = L^i G K_{n-2} $ and so on.  \Qed

 \medskip
\subsection{The main induction}
 
Recall that the set 
$\seq_{n}^{r}\{1,2 \}$   
is a basis of a permutation submodule of $V^{\otimes n}$, 
which we now denote $ \permmod{2r-n}{n}$ 
(the argument $2r-n$ counts the excess of $1$'s over $2$'s). 
Now let 
$$ v(r,n) \; := \;  111 \ldots 11 222 \ldots 22 \; \in \seq_{n}^{r}\{1,2 \} $$
Then $ v(r,n) $ generates $ \permmod{2r-n}{n}$ as a $ T_{n}$--module.


Our argument for showing that $ V^{\otimes n } $ is tilting will be 
an induction on $n$. 
The inductive step is based on the following Lemma:
\begin{lem}
Assume that $n\geq 2$. 
Then there are isomorphisms in $T_{n-2} \, \mbox{--mod} $:
\eql(indstep1) 
F (V^{\otimes n } ) \cong V^{\otimes n-2 }  
\eq
\eql(indstep2) 
F\permmod{s}{n} \cong  \left\{ \begin{array}{ll} 
\permmod{s}{n-2}   & |s| < n
\\
0 & |s|=n
\end{array} \right. 
\eq
\end{lem}
{\em Proof:} $F$ is multiplication by $ \idem
= \frac{1\,\,}{[2]_q} U_{n-1} $, which acts 
in the last two factors of $ V^{\otimes n } $ through the matrix:
$$ \frac{1\,\,}{[2]_q}\mat{cccc} 
0&0&0&0\\ 
0&q&1&0 \\
0&1&q^{-1}&0 \\
0&0&0&0 \tam$$
which has eigenvalues $1$ with multiplicity 1 and $0$ with multiplicity 3 
and therefore is a projection onto a one dimensional space. 
Let $ w_2 \in V^{\otimes 2 } $ be an eigenvector to eigenvalue $1$
(say $ w_2 = q \, v_1 \otimes v_2 + v_2 \otimes v_1 $). 
Then, 
cf. proposition~\ref{Vres},  
the (inverse of the) first   
isomorphism is given by 
$ v_{i_1} \otimes  v_{i_2} \otimes \ldots \otimes  v_{i_{n-2}}  \mapsto
 v_{i_1} \otimes  v_{i_2} \otimes \ldots \otimes   v_{i_{n-2}} \otimes w_2 $.

To get the 
second
 isomorphism, note first that the above map clearly induces
$$ F\permmod{s}{n} \subseteq  \permmod{s}{n-2}$$
But then equality follows from the 
first     
isomorphism combined with 
$ V^{\otimes n } = \bigoplus_s \permmod{s}{n} $ 
and the analogous formula for
$ V^{\otimes n-2 } $.
\Qed

\medskip

Note that this Lemma relates level $n$ to $n-2$. 
Accordingly our inductive argument
for proving that $V^{\otimes n}$ is tilting will require {\em two} base cases: 
$ n=1 $ and $ n=2$. 
Both are straightforward 
($n=1$ it trivial, and for $n=2$ we have 
$ V^{\otimes 2} = 3 \Delta_2(2) \oplus \Delta_2(0)  
=
3 \nabla_2(2) \oplus \nabla_2(0)  $   
by the proof of the Lemma). 


We now consider the exact sequence 
\eql(four terms)
 0 \rightarrow K_n \rightarrow  G \circ F( V^{\otimes n}) 
\stackrel{\varphi_n}{\rightarrow}
V^{\otimes n} \rightarrow 
C_n \rightarrow 0 
\eq
where $ \varphi_n $ is the adjointness map. 
\prl(inductive step)
Assume that $ \varphi_n $ is injective for all $ n $. Then 
$ V^{\otimes n } $ and $ \permmod{r}{n}$
are tilting modules for $ T_n $ for all $n,r $.
\end{pr}
{\em Proof}. Since $ F \circ G = Id $, we have 
that $ F \circ G \circ F = F $ and hence $ F(C_n) = 0 $. 
Thus, cf. Proposition~\ref{Apply F}, $C_n$ has 
only (copies of) the `trivial' one
dimensional module (= $ \Delta_n(n) $) as composition factors.  
But then 
$ C_n $ is semisimple by quasiheredity (or otherwise). 

Now work by induction on $n$. 
By the Lemma and the inductive hypothesis 
$F (V^{\otimes n } )$ is tilting, so 
in particular  $F (V^{\otimes n } )$ has a standard filtration, and
then so does $G \circ F (V^{\otimes n } )$ by Proposition~\ref{apply G}.
But $ \varphi_n $ is assumed to be injective, 
so $K_n = 0$ and (\ref{four terms})
becomes a short exact sequence with $ V^{\otimes n } $ 
in the middle and with extremal terms in $ {\cal F }_n(\Delta) $. 
But then $ V^{\otimes n } $ too lies in $ {\cal F }_n(\Delta) $.
Since the matrices 
representing the action of $ T_n $ on $ V^{\otimes n } $ are selfadjoint with 
respect to the canonical, non--degenerate form 
(note from the presentation that the algebra is isomorphic to its opposite), 
$ V^{\otimes n } $
is contravariant selfdual and so  tilting. 
But then also $ \permmod{r}{n} $ is tilting
as a direct summand of $V^{\otimes n }$
\Qed
\medskip



We have that $\varphi_n$ injective implies $V^{\otimes n}$ tilting.  
As an aside we note that the reverse implication also holds, 
and rather more generally. 
Indeed for  $ M \in {\cal F }(\Delta) $, the adjointness
map $ G \circ F(M) \rightarrow M $ is injective. This is clear if 
$ M \cong  \Delta $ and otherwise it follows by induction on the number of 
$ \Delta $-factors in $M$ using the following commutative diagram
$$
\begin{array}{cccccccccc}
 &  & 0 &  \rightarrow &  0 & \rightarrow 
&  0 &  &  
 &    \\

 &  & \downarrow  &   &  \downarrow  & 
&  \downarrow  &  &  
 &    \\

0 &  \rightarrow & G \circ F (C) &  \rightarrow & G \circ F(M) &  
\rightarrow
&  G \circ F (\Delta) &  \rightarrow &  
 &  0  \\

 &  & \downarrow  &   &  \downarrow  & 
&  \downarrow  &  &  
 &    \\

0 &  \rightarrow &  C &  \rightarrow & M &  
\rightarrow
&   \Delta &  \rightarrow &  
 &  0  
\end{array}
$$ 
where standard $\Delta$ is such that $M \twoheadrightarrow \Delta$,
noting that the second row of the diagram is exact because of 
equation~(\ref{dtool}). 


\medskip

Now  spelling out the definitions, $ \varphi_n $ is the multiplication map 
$$ \varphi_n: T_n U_{n-1}
 \otimes_{{U_{n-1}T_n U_{n-1}}}  U_{n-1}
V^{\otimes n } 
 \rightarrow  V^{\otimes n }  $$

\medskip
The rest of 
the construction of our inductive step amounts to 
a careful combinatorial analysis of this map. 
First of all consider $ T_n U_{n-1}$ as a right module over
$ {{U_{n-1}T_n U_{n-1}}} $. 
As such it is easy to see that it is generated by the elements 
\eql(generators) 
 U_{n-1},\,\;\;\;\;   U_{n-2} \, U_{n-1}, \,\;\;\;\;
\ldots ,\,\;\;\;\; U_{1}\,\cdots  \, U_{n-2} \,U_{n-1} 
\eq
But then any element of $ T_n U_{n-1}
 \otimes_{{U_{n-1}T_n U_{n-1}}}  U_{n-1}
V^{\otimes n } $ can be represented in the form
\eql(bythis) 
\sum_{k} \, U_{k}\,\cdots  \, U_{n-2} \,U_{n-1} 
\otimes_{{U_{n-1}T_n U_{n-1}}}  v_k 
\eq
for some $ v_k \in U_{n-1}  V^{\otimes n }$. 
We must show that this expression is 
zero if its image under the multiplication map is zero, i.e. if 
$$ \sum_{k} \, U_{k}\,\cdots  \, U_{n-2} \,U_{n-1} \, v_k = 0  $$
To do this, the following notation will be useful:
\newline 
Recall that 
$ i_1 i_2 \ldots i_n  \in \seq_{n}\{1,2 \}$   
is a basis element of $ V^{\otimes n } $. 
Denote by 
$ i_1 i_2 \ldots i_{k-1} \underline{12}\, i_{k+2}   \ldots i_n $ 
the vector 
$ U_k\, ( i_1 i_2 \ldots i_{k-1} 12\, i_{k+2} \ldots i_n)  $. 
In other words: 
$$ 
i_1 i_2 \ldots i_{k-1} \underline{12}\, i_{k+2}  
\ldots i_n 
=
q\, v_{i_1} \otimes v_{i_2} \otimes \ldots 
\otimes v_1 \otimes v_2 \otimes \ldots \otimes v_{i_n} +
    v_{i_1} \otimes v_{i_2} \otimes \ldots
\otimes v_2 \otimes v_1 \otimes \ldots \otimes v_{i_n} 
$$

\medskip
To establish usage of this notation 
we first work out a couple of low rank examples. 

\noindent
{\bf Example 1.} Consider $n=3$.  
Then $\{ U_1, U_2 \} $ generate $T_n$ and
$ V^{\otimes n} $ has dimension 8 and is the direct sum of 4 permutation 
modules. The injectivity of our map can be checked on each of them.
On the one dimensional permutation modules the statement is trivial since
$ G \circ F $ kills them.
Let us then consider the permutation module generated by $112 $ 
(that generated by $122$ is isomorphic to it). 
Applying 
$ G \circ F $, we get by (\ref{bythis}) the two vectors 
$$ 
U_2 \otimes_{U_2 T_3 U_2} 112 
= 1 \otimes_{U_2 T_3 U_2} U_2 112
= 1 \otimes_{U_2 T_3 U_2} 1 \underline{12} , \hspace{.5cm}
U_1U_2 \otimes_{U_2 T_3 U_2} 112 
=
    U_1  \otimes_{U_2 T_3 U_2} 1 \underline{12}   
$$
(NB, $U_2 \otimes 211 = 0$ and 
$U_2 \otimes 121 \propto U_2 \otimes 112$). 
The images under the multiplication are 
$$  
1 \underline{12} , \,\,\,\,
    \underline{12}\, 1 ,   \,\,\,\    
$$
and these are linearly independent, so also in this case the 
statement is clear.
\medskip

\noindent
{\bf Example 2}: Consider $ n=4$. Then $T_n$ 
is generated by $\{U_1, U_2 , U_3 \}$ and $ V^{\otimes n } $ is the 
sum of five permutation modules.

\medskip
Let us consider the permutation module$\, \permmod{0}{4} $ 
generated by $ 1122 \in V^{\otimes 4 } $. Using (\ref{bythis}),
$ G \circ F  \permmod{0}{4}$ is spanned by the vectors
$$ 
\begin{array}{ll}
1 \otimes 12 \underline{12} , &  
1 \otimes  21 \underline{12} ,
\\
U_2 \otimes 12 \underline{12} , & 
U_2 \otimes 21 \underline{12} , 
\\
U_1 \,  U_2 \otimes 12 \underline{12} , &
U_1 \,  U_2 \otimes 21 \underline{12} 
\end{array}
$$
i.e. by the set $\{ U_1 U_2 \otimes w \underline{12} , \; 
U_2 \otimes w \underline{12} , \; 
1  \otimes w \underline{12}  \; 
\; | \; w \in \seq_{2}^{1}\{1,2 \}  \}$. 

The images under the multiplication map are the vectors 
$$ 
\begin{array}{ll}
 12 \underline{12} , &  
 21 \underline{12} , 
\\
 1 \underline{12} 2 , & 
 2 \underline{12} 1 , 
\\
   \underline{12} 12 , &
   \underline{12} 21 . 
\end{array} 
$$
Note that these vectors are not independent, since 
$$ 
q 12 \underline{12} + 21 \underline{12} 
=\underline{1 2}\, \underline{1 2} = q \underline{12} 12 +
\underline{12} 21
$$ 
On the other hand this {\it `trivial'} dependency is the only one, 
as can seen by a dimension counting: our permutation module has dimension 
$\left( \begin{array}{c}
4 \\ 2 \end{array} \right) = 6$ 
and is generated by the above vectors together with $ 1122 $. So 
there is exactly one relation between them, the one we have pointed out.

Since there is a corresponding dependency amongst the first set of vectors:
\[
q1 \otimes 12 \underline{12} + 1 \otimes 21 \underline{12}
= 1 \otimes  \underline{12}\underline{12} = 
U_1 U_2 \otimes  \underline{12}\underline{12} = 
qU_1 U_2 \otimes 12 \underline{12} + U_1 U_2 \otimes 21 \underline{12}
\]
our claim is proved in this case as well.
\medskip

We now turn to the general case. 
\begin{theo} 
$ V^{\otimes n } $ and
$\permmod{r}{n} $ are 
tilting modules for all $ n $.
\end{theo}
{\em Proof.} 
By proposition~\ref{inductive step} it is enough to show that 
$\varphi_n $ is injective for all $ n$. 
$ F $ and $ G $ are additive functors, so the claim 
can be verified on the permutation submodules. 
By (\ref{bythis}), $G \circ F$ on 
$\permmod{2r-n}{n} $
is spanned by 
$$
\{ \; X \otimes w  \underline{12} \;\; 
         | \;\; X \in \{ U_1 .. U_{n-1},\; U_2 .. U_{n-1},\; .., \;U_{n-1},\; 1 \} 
         ; \;\;\; w \in \seq_{n-2}^{r-1}\{1,2 \} \}  
$$


Let us denote by $\SSv^r$ the set of vectors of the form 
$$  
i_1 \, i_2 \, i_3 \,  \ldots \underline{12} \, \ldots  
i_{n-2} \, i_{n-1} \, i_n  \mbox{ \,\,\,\, } i_k \in \{ 1,2 \} 
$$  
inside $\permmod{2r-n}{n} $.
Note that these are the images of the above vectors under the multiplication
map. 
A simple counting argument shows that 
$ |\SSv^r |   =  
\left( \begin{array}{c}
n-1 \\ 1 
\end{array} \right)
\left( \begin{array}{c}
n-2 \\ r-1
\end{array} \right)$. 
Note that $\permmod{2r-n}{n} $
is spanned by $\SSv^r \cup \{ v(r,n) \}$. 
But $\SSv^r$ is not linearly independent. 
To each vector in $\permmod{2r-n}{n} $
of the form 
$$  
i_1 \, i_2 \, i_3 \,  \ldots \underline{12} \, \ldots  \underline{12} 
\ldots \, i_{n-2} \, i_{n-1} \, i_n 
$$
we may associate a dependency:
\eql(bythat)
q \ldots 12 \ldots \underline{12} \ldots +\ldots 21 \ldots \underline{12} \ldots
=\ldots \underline{12} \, \ldots  \underline{12} \ldots 
=q\ldots \underline{12} \ldots 12 \ldots+\ldots \underline{12}\ldots 21 \ldots
\eq 
Note that each of these dependencies has a preimage in 
$G\circ F \permmod{2r-n}{n} $:
\[ 
q U_j..U_{n-2}\otimes \ldots 12 \ldots w_j w_{j+1} \ldots \underline{12}
\; + U_j..U_{n-2}\otimes \ldots 21 \ldots w_j w_{j+1} \ldots
\underline{12}
\] \[
=U_j..U_{n-2}\otimes \ldots\underline{12}\ldots w_j w_{j+1} \ldots\underline{12}
\] \[
=U_i..U_{n-2}\otimes \ldots w_i w_{i+1}\ldots\underline{12} \ldots\underline{12}
\] \[
=qU_i..U_{n-2}\otimes \ldots w_i w_{i+1} \ldots 12 \ldots\underline{12}
\; +U_i..U_{n-2}\otimes \ldots w_i w_{i+1} \ldots 21
 \ldots\underline{12}  
\] 
(recall that we are tensoring over $ \idem T_n \idem
 = U_{n-1} T_n U_{n-1} $).


Let $\SSv'$ denote the subset of $\SSv^r$ 
in which no subsequence 12 appears before the $\underline{12}$. 
We will now show, using the set of linear dependencies above, 
that all vectors not in the subset $\SSv'$ may be
discarded without affecting the spanning property, 
i.e. $\permmod{2r-n}{n}$ is also spanned by $\SSv' \cup \{ v(r,n) \}$. 

Let the `underlying' sequence $ u(s) \in \seq_n\{1,2 \} $ of 
$ s \in \SSv^r $ be the sequence obtained by removing the underline from $s$.
Note that any $ s \in \SSv^r $ of the form 
$..12..\underline{12}..$ can be written as a linear combination
of $..\underline{12}..12..$ and  $..21..\underline{12}..$ 
and $..\underline{12}..21..$ using (\ref{bythat}), 
and that  the last two have underlying sequences 
later in the \lex\ 
order than $..12..\underline{12}..$. 
Let $ s = ..12..\underline{12}.. \in \SSv^r \setminus \SSv' $. 
There may be other pairs 12 before $\underline{12}$  in $s$, but 
we may take it that the  pair 12 written explicitly is the 
{\em  first} such. 
Using  the above remark we then replace $s$ by a 
linear combination of an element of $ \SSv' $ and elements
of $\SSv^r$ not necessarily in $ \SSv' $ but whose underlying sequence is 
later in the \lex\ 
order than $u(s)$. 
Iterating, we arrive at $ 22.21\underline{12}11..1 $ 
(or similar) which is in  $ \SSv' $.

Note that $\SSv'$ has a natural bijection with 
$\seq^r_n \{1,2 \} \setminus \{ 22..211..1 \}$ 
(every element of this set has at least one subsequence 12 
--- just underline the first of these). 
Thus  $\SSv' \cup  \{ v(r,n) \}$  is a basis
of $\permmod{2r-n}{n}$, and in particular it is linearly independent. 
This means that all linear dependencies in $\SSv$ can be constructed 
from those of form (\ref{bythat}). 
But each of these dependencies has a preimage in 
$G\circ F \permmod{2r-n}{n}$, so $\varphi_n$ has a trivial kernel. 
\Qed


\begin{cor} $ V^{\otimes n } $ is a full tilting module for $ T_n $. 
\end{cor}
{\em Proof.}
We proved in the Theorem that 
$ \permmod{s}{n} $ is a tilting module for all $ n,s $.
Now we have the restriction rule 
$$ 
\mbox{Res}^{T_n}_{T_{n-1}} \permmod{s}{n} = \permmod{s-1}{n-1} \oplus 
\permmod{s+1}{n-1} \,\,\, 
\mbox{ for } s \in \{- n,-n+2 \ldots, n-2,n \} 
$$
Combining this with the restriction rule for the standard modules 
\cite{Martin91} 
$$ 
[\mbox{Res}^n_{n-1}\Delta_n(s)]= 
[\Delta_{n-1}(s+1)] +[ \Delta_{n-1}(s-1)]
\hspace{.2cm} \mbox{ for } s \in \{0/1\ldots, n-2,n \} 
$$
it is easily proved by induction that 
$$
[\permmod{s}{n}] = [\Delta_{n}(s)]+[\Delta_{n}(s+2)]+ \ldots 
\hspace{.2cm}  (s \geq 0) 
$$
In other words $ (\permmod{s}{n}: \Delta_{n}(u)) = 1 $ for $u$ less than 
or equal $ s $ in the heredity order; $0$ otherwise. 
But then the tilting module $ \tilting{s}{n} $ must occur as a summand of
$ \permmod{s}{n} $ (with multiplicity one). \Qed



\section{The blob crypto--tensor module case}\label{blobtilt}

\medskip 
Let us now consider the blob algebra situation. 
Our overall strategy will be  similar to the one used in
the previous section. 
In particular the quasiheredity arguments carry over almost unchanged. 
On the other hand, the
actual calculation requires some new combinatorial ideas. 

\medskip
We keep the condition that $ [2]_q \not= 0 $, but assume also that 
$[m]_q \not=0 $ (where $ m $ is as in  section~\ref{Preliminaries}). 
In that case the blob algebra $ b_n^k = b_n $ is \qh, 
see \cite{MartinWoodcock2000}. 
In fact, setting $\idem = \frac{1}{[2]_q} U_{n-1} $, we have as for $T_n(q)$ that 
$ \idem $ is part of a heredity chain for $b_n $ 
and   
$ \idem b_n \idem \cong  b_{n-2} $.
Accordingly the results from the 
previous section involving quasiheredity hold in this case as well.
We state them here indicating the necessary modifications of the
previous proofs.


As before we have a localization functor $ F $ 
$$ 
F: b_n \, \mbox{--mod} \rightarrow b_{n-2} \, 
\mbox{--mod}:\,\, M \mapsto \idem M 
$$ 
and a  
globalization functor $G$
$$ 
G: b_{n-2} \,\mbox{--mod} \rightarrow b_{n} \,\mbox{--mod}:\,\, 
N \mapsto b_{n} \idem
\otimes_{ b_{n-2} } N $$

We denote as before 
the standard (costandard etc.) modules in $ b_n \mbox{-mod} $ by 
$ \Delta_n(\lambda) $ ($\nabla_n(\lambda)$ etc.), but in this case the 
parametrizing set is $ \Gamma_n =\{-n, -n+2, \ldots, n-2, n \} $ 
\cite{MartinSaleur94a,MartinWoodcock2000}. 
We then have the following version of Proposition~\ref{Apply F}
\prl(blobApply F)
Assume that $ \lambda \in \Gamma_{n-2} $. Then 

\noindent
i) $F L_n(\lambda) = L_{n-2}(\lambda) $. 

\noindent
ii) $F \Delta_n(\lambda) = \Delta_{n-2}(\lambda) $ and 
$F \,  \nabla_n(\lambda) = \nabla_{n-2}(\lambda) $.

\noindent
iii) $ F P_n(\lambda ) = P_{n-2}(\lambda )$ and $ F I_n(\lambda ) 
= I_{n-2}(\lambda ) $.

\noindent
Otherwise ( i.e. for $ \lambda = n \mbox{ or }\lambda = -n $ ) we have $F L_n(\lambda) = 
F \Delta_n(\lambda) = F \nabla(\lambda ) = 0$.
\end{pr}
We then get as before
\prl(blobapply G) 
Supposing a module $ M \in {\cal F }_{n-2}(\Delta)\, $ 
(i.e. $\, M \in b_{n-2} \,\mbox{-mod}$ has a standard filtration), 
then $G(M)$ also has a standard filtration. 
Furthermore 
$$
( G(M) : \Delta_{n}(\mu) ) = \left\{
\begin{array}{ll}
( M : \Delta_{n-2}(\mu) ) &  \mbox{if  }  \mu \in \Gamma_{n-2} \\
0 & \mbox{otherwise}  \\
\end{array}
\right.
$$
\end{pr}
{\em Proof.} The proof is once again an application of the cohomological 
criterion for standard filtrations. One shows that
$ \mbox{ Ext}_{T_{n}}^1 ( G(M),\nabla_n(\mu ) ) = 0 $ for all $ \mu $. The 
special case $ \mu = n $ from the Temperley-Lieb situation now becomes two special 
cases $ \mu = n $ and $ \mu =- n $, each of which can be treated as before. \Qed

The cohomological statement (\ref{dtool}) also carries over:

\begin{co}
\eql(blobdtool) 
L^{i}G \Delta_{n-2}(\lambda)  \,\, = \,\,
\left\{
\begin{array}{ll}
\Delta_n(\lambda) & \mbox{    if } i = 0 \\
0      & \mbox{    otherwise }
\end{array}
\right.
\eq 
\end{co}

The set $\seq_{2n}^r\{1,2\}$ is a basis of a permutation module
of $\ourrep$ which we denote $\permmod{2r-2n}{n}$. 
(For example a basis of $\permmod{0}{2}$ is 
$\{  1122, 1212, 1221, 2112, 2121, 2211 \}$.) 
Evidently  
\beqa \label{permres1}
\Res(b_{n},b_{n+1},{\permmod{\lambda}{n+1}})
 &=& \permmod{\lambda +2}{n} \oplus 2 \permmod{\lambda}{n} 
     \oplus \permmod{\lambda -2}{n}
\eeqa
We also have a blob algebra version of (\ref{indstep1}): 
\begin{lem} \label{lemma 2}
Assume that $n\geq 2$. Then there 
are isomorphisms in $b_{n-2} \, \mbox{--mod} $:
\eql(blobindstep1) 
F (V^{\otimes 2n } ) \cong V^{\otimes 2( n-2) }  
\eq
\eql(blobindstep2) 
F\permmod{s}{n} \cong  \left\{ \begin{array}{ll} 
\permmod{s}{n-2}   & |s| < 2n-2
\\
0 & |s|=2n, \; 2n-2
\end{array} \right. 
\eq
\end{lem}
{\em Proof.} The Temperley--Lieb argument goes through almost unchanged:
$ F $ is multiplication by the idempotent $ \idem = \frac{1}{[2]_q} U_{n-1} $ 
which acts only in the first two and last two factors of $ V^{\otimes 2n } $. 
The isomorphism $ V^{\otimes 2( n-2)} \rightarrow F (V^{\otimes 2n } ) $ is given
by $ w \in   V^{\otimes 2( n-2)} \mapsto  ev_1 \otimes  w \otimes ev_2 $ 
where $ ev_1 \otimes ev_2 \in V^{\otimes 4} $ is an eigenvector to eigenvalue $1$
of our idempotent $ \idem $. \Qed

\medskip
\noindent
We consider also in the blob algebra setting the adjointness map 
$ \varphi_n:G \circ F( V^{\otimes 2n}) \rightarrow V^{\otimes 2n} $ 
and 
get a four term exact sequence:
\eql(blobfour terms)
0 \rightarrow K_n \rightarrow  G \circ F( V^{\otimes 2n}) 
\stackrel{\varphi_n}{\rightarrow}
V^{\otimes 2n} \rightarrow 
C_n \rightarrow 0 
\eq

\prl(blobinductive step)
Assume that $ \varphi_n $ is injective for all $ n $. Then 
$ V^{\otimes 2n } $ and $ \permmod{r}{n}$
are tilting modules for $ b_n $ for all $n,r $.
\end{pr}
{\em Proof.} 
The Temperley--Lieb proof carries over almost verbatim. \Qed

\medskip
We are then once again left with the task of showing injectivity of 
the adjointness map 
$ \varphi_n: G \circ F (V^{\otimes 2n }) \rightarrow V^{\otimes 2n }  $. 
Let us as before first work out a low rank example 
(which will this time also be needed in the general argument).

\medskip
\noindent
{\bf Example 1.} Consider $n=2$. Then 
the blob algebra $ b_2 $ acts on $ V^{\otimes 4 } $ as described in 
section~\ref{the rep}. 
For example
\[
U_1 1212 
\;\; = \;\;
st 1212 + s 1221 + t 2112 + 2121
\;\;\; =: \;\;  \underline{12}\, \underline{12}
\]
(NB, the underline notation is refined here to accommodate the
definition of $\ourrep$ --- 
the {\em position} of the underline 
to left or right of centre determines precisely which linear
combination it corresponds to). 
A generating set of $ b_2 $ is $ U_1 $ and $ U_0 =a^2 \egen$. 
Hence the vectors 
$$   1  \otimes  \underline{12}\, \underline{12}, \,\,\,
     U_0  \otimes  \underline{12}\, \underline{12} $$
generate $ G \circ F  ( V^{\otimes 4 }) $.
The images under the multiplication map $\varphi_2$ are 
$$    
\underline{12}\, \underline{12}, \,\,\,\,\,\,\,\, 
\underbrace{ 1   \underline{12} 2} 
$$ 
(where $ \underbrace{ 1   \underline{12} 2}$ denotes a certain linear
combination of $  1   \underline{12} 2 $ and $ 2  \underline{12} 1$:
$U_0 \underline{12}\, \underline{12} 
= r 2121 +2211 +rst 1122 +st1212$)
which are independent.

\medskip
More generally, a spanning set for $ G \circ F ( V^{\otimes 2n }) $ is
given by
$$
B_n \; := \; \{ X \otimes \underline{12} w \underline{12}
\;\; | \;\; X \in 
\left\{ \begin{array}{c} 
1, \\ U_{n-2}, \\ U_{n-3}U_{n-2}, \\ \vdots \\ U_0 U_1 .. U_{n-2}
\end{array} \right\} 
;
w \in \seq_{2n-4}\{1,2\} \}
$$
The images under the multiplication map are vectors of the form 
$ .. \underline{12} .. \underline{12} .. $  where the concatenation
of the subsequences indicated by ellipsis is the sequence $w$ from 
$\seq_{2n-4}\{1,2 \}$, and the first and second 
$\underline{12}$ are equidistant
from the left and right hand end respectively; and 
$ .. \underbrace{ 1   \underline{12} 2} ..$, where the concatenation
of the subsequences indicated by ellipsis is the sequence from 
$\seq_{2n-4}\{1,2 \}$. 
This is a consequence of the following straightforward exercise in 
the blob algebra relations:
\begin{lem}
Let $ \idem= \frac{1}{[2]_q} U_{n-1} $. 
The set $ \{ 1,\,\; U_{n-2},\,\; U_{n-2}U_{n-3},\, \ldots , \,
U_{n-2}U_{n-3} \ldots U_0 \} $ generates $ b_n \idem $ as a right 
$ \idem b_n \idem $--module.
\end{lem}

It will be evident that linear dependencies arise in general between
these vectors, in a way analogous to the $T_n$ case. For
example it is easy to write down a linear dependence involving 
$  \underline{12} 1212 \underline{12}$ and 
$ 12 \underline{12} \underline{12} 12$ (and others).

It is straightforward, using this machinery, to verify that 
$\varphi_n$ is injective and hence
our module is tilting, up to $b_4$ and a little beyond. 
We will prove injectivity for general $n$ by a slightly different route. 

Define first numbers $ v(n) $ by the recursion 
$v(0)=1, v(1)=1, v(-1)=3 $ and
$$ 
v(n)= 4 v(n-1)-v(n-2) \;\;\;\; \mbox{  if } n \geq 2,
\,\,\,\,\;\;\;\; 
v(n)= 4 v(n+1)-v(n+2) \;\;\;\; \mbox{  if } n \leq -2 
$$
\prl(rank1) Let $ \rho_{}(n) $ be the representation $\ourrep$ of $ b_n $ on 
$ V^{\otimes 2n } $. Then

\noindent
1) $ \varphi_n: G \circ F(\rho_{}(n)) \rightarrow \rho_{}(n) $ is
injective, so $ \rho_{}(n) $ is a tilting module. 

\noindent
2) $ ( \rho_{}(n) : \Delta_n(\lambda) ) = v(\lambda) $.

\noindent
3) Set $r_n = \dim (G \circ F (\rho_{}(n)) )$. Then $r_1 =0$, 
$r_2 =2$, and 
\[ r_n = 4 r_{n-1} + 4^{n-2} - r_{n-2} . \] 

\end{pr}
{\em Proof:}
By induction on $n$. The case $ n=1$ follows easily from the fact that 
$$ 
\rho_{}(1) = V \otimes V \cong \Delta_1(1) \oplus 3\Delta_1(-1) 
$$ 
and {\em 1)} and {\em 3)} 
of the $ n= 2 $ case is the calculation done above. 
The calculation also shows that $ \dim F(\rho_{}(2)) = 1 $, so we get that
$ \Delta_2(0) $ occurs in $\rho_{}(2)$ with multiplicity $1$. 
We can then read off the other two multiplicities using 
$ \res^{b_2}_{b_1} \, \rho_{}(2) = 4 \rho_{}(1) $ 
and the restriction rules 
for the standard modules, thus verifying {\em 2)}. 
For the reader's convenience we express this last point in formulas: 
write $\rho_{}(2)$ in the Grothendieck group as follows
$$ 
\rho_{}(2) = a_2 \,\Delta_2(2) +a_0\, \Delta_2(0) + \, 
a_{-2} \,\Delta_2(-2) 
$$ 
Applying $ F $ to this expression we get $ a_0=1$ and applying the
restriction functor to it we get 
$$ 
4\, \rho_{}(1) = \res \rho_{}(2) = (a_2 + a_0 )\, \Delta_1(1) + 
(a_{-2} + a_0 )\, \Delta_1(-1) 
$$
and {\em 2)} now follows from 
$ \rho_{}(1) =  \, \Delta_1(1) +3 \, \Delta_1(-1) $.

Now assume the Proposition for $ n^{\prime}$ with $ n^{\prime} < n $. 
Then 
$ F( \rho_{}(n)) = \rho_{}(n-2) $ is tilting and 
$$ 
\SCh{ F( \rho_{}(n)) : \Delta_{n-2}(\lambda) } \;\; = \;\; 
  \SCh{ \rho_{}(n-2) : \Delta_{n-2}(\lambda) }  = v(\lambda) \,\,\, \mbox{ if }
|\lambda| \leq n-2 
$$ 
But then also 
$$ 
\SCh{ G \circ F( \rho_{}(n)) : \Delta_{n}(\lambda) } 
  = v(\lambda) \,\,\, \mbox{ if } |\lambda| \leq n-2 
$$ 
since $ G $ is exact on $ {\cal F}(\Delta) $ and takes standard modules 
to standard modules.
Note \cite{MartinWoodcock2000} the short exact sequence
\beqa \label{S res}
0 \rightarrow \Delta_{n-1}(\lambda \pm 1)  
  \rightarrow \mbox{Res}^{b_{n}}_{b_{n-1}} \Delta_{n}(\lambda)  
  \rightarrow \Delta_{n-1}(\lambda \mp 1) \rightarrow 0
\hspace{.5cm} (\lambda \mp 1 
\stackrel{\hlt}{\mbox{{\scriptsize $\hgt$}}} 
\lambda)
\eeqa
($\Delta_{n-1}(\nu)$ to be interpreted as 0 if $|\nu| > n-1$). 
We can then calculate $ r_n $ as follows
$$ 
r_n= \sum_{\lambda: |\lambda| \leq n-2} \, v(\lambda) 
| \Delta_{n}(\lambda) | = 
 \sum_{\lambda: |\lambda| \leq n-2}  v(\lambda) 
\left(\,| \Delta_{n-1}(\lambda+1) | + | \Delta_{n-1}(\lambda-1) |\, \right) =$$
$$ \left( v(n-2) +v( n-4)\right) | \Delta_{n-1}(n-3) | +  \ldots +
 \left( v(-n+2) +v(- n+4)\right) | \Delta_{n-1}(-n+3) | $$
$$ +v(n-2) | \Delta_{n-1}(n-1) | + v(-n+2) | \Delta_{n-1}(-n+1) |  $$ 
Using the recursion formula for $v(n) $, this becomes
$$ 
4 \left\{ v(n-3) | \Delta_{n-1}(n-3) | + \ldots 
      +  v(-n+3) | \Delta_{n-1}(-n+3) | \right\} + v(n-2) + v(-n+2) 
$$

We now apply the induction hypothesis (part {\em 2)} and {\em 3)}) 
and get that the first term is equal to
$ 4 r_{n-1} $, while the sum $ v(n-2) + v(-n+2)$ is equal to
$$ 
\dim \rho_{}(n-2)/ G \circ F(\rho_{}(n-2)) = 4^{n-2} - r_{n-2} 
$$
Combining this we have shown {\em 3)} at level $n$.

To prove {\em 2)} at level $n$ note first that $ F $ annihilates   
$ \rho_{}(n)/ G \circ F(\rho_{}(n)) $, so 
it can be written in the Grothendieck group as a sum
$$ 
a_n \,\Delta_n( n) + a_{-n}\, \Delta_n(-n) 
$$
We restrict and apply the formula 
$ \res \,\rho_{}(n) = 4 \, \rho_{}(n-1) $, and find by comparing coefficients 
that
$$ 
a_{n}+v(n-2) = 4 \, v(n-1) \;\;\; \mbox{ and   } \;\;\; 
 a_{-n}+v(-n+2) = 4 \, v(-n+1)
$$
Since the $\Delta $--multiplicities of  
$G \circ F(\rho_{}(n)) $ are already known by induction, this 
shows {\em 2)} at level $n$.


Let us now finally prove {\em 1)} at level $n$. 
This is the most tricky part of our proof 
and involves some interesting combinatorics on sequences. 
Let $ v \in \varphi_n(B_n) $, then 
the four sequences which occur as summands of $v$ are either of the
form 
\[
 \{ x 12 y 12 z, \; x 12 y 21 z, \; x 21 y 12 z, \; x 21 y 21 z \}
\]
where $x,y,z \in \seq_{}\{1,2 \}$, or 
\[
\{ x 11  22 z , \;  x 12  12 z , \; x 21  21 z , \; x 22  11 z \}
\]
where $x,z \in \seq_{}\{1,2 \}$. 
Define $ u(v) \in \seq_{2n}\{1,2 \} $ to be the 
\lex ally earliest 
sequence that occurs as a summand of $ v $. 

One easily sees from the description of the elements of 
$\varphi_n( B_n ) $ that $ u(v) $ 
satisfies the rule 
\eql(* rule)  
    au(v)b = u(avb) 
\eq 
for $ a,b  \in
               \{1,2 \} $ such that $ avb \in \varphi_n( B_{n+1} ) $.


Now define 
for all $ n \geq 1 $ a subset $ E_n $ of our 
representation space $ V^{\otimes 2n } $ as follows:
$$ 
E_1 := \emptyset, \,\,\,\;\;\;\; E_2:=\{ \underline{12} \underline{12},
\underbrace{ 1   \underline{12} 2} \} 
$$
then for $ n \geq 2 $ : 
$$ 
E_n^1:= \{ \,1x1,  1x2, 2x1, 2x2 \, |\,  x \in E_{n-1}\, \} 
$$
$$ 
E_n^2:= \{ \, \underline{12} w \underline{12} \, |\, w \in 
\seq_{2n-4}\{1,2 \} \setminus  u(E_{n-2}) \,  \} 
$$
$$ 
E_n:=E_n^1 \cup  E_n^2 $$
Consider now the following properties of $E_n$:

\noindent
{\bf Claim:} i)  $| E_n| = r_n$ 

\noindent
ii) $E_n$ is a basis of $ \varphi_n( G \circ F(\rho_{}(n))) $

\noindent
iii) $|E_n| = | u(E_n )| $

\noindent 
iv)  $u(E_n^1 ) \cap  u(E_n^2 ) = \emptyset$

\medskip
Part {\em 1)} of the Proposition is a consequence of i) and ii) 
since we already know that  $ \dim  G \circ F(\rho_{}(n)) = r_n $. 
In order to prove the claim we again proceed by induction. 
Since $ E_n^i$ are only defined for $ n \geq 3 $
we take $ n=3 $ as base of the induction, but actually 
i), ii) and iii) also make sense for $n=1, 2 $ and 
basically follow from the calculations
prior to the Proposition: note that 
$$   
u( \underline{12}\, \underline{12})= 1212, \,\,\,\,\,\,\,\, 
u( \underbrace{ 1 \underline{12} 2} )= 1122 
$$ 
to obtain iii). Now for $ n=3 $ we have 
$$ 
E_3^1 = \left\{ 
\begin{array}{cccc}
1 \underline{12}\, \underline{12}1, \\
1 \underline{12}\, \underline{12}2, \\
2\underline{12}\, \underline{12}1, \\
2\underline{12}\, \underline{12}2, 
\end{array}
\begin{array}{cccc}
1\underbrace{ 1 \underline{12} 2}1, \\
1\underbrace{ 1 \underline{12} 2}2, \\
2\underbrace{ 1 \underline{12} 2}1, \\
2\underbrace{ 1 \underline{12} 2}2 
\end{array}
\right\}
$$
while $$ E_3^2 =
 \left\{ 
\begin{array}{cccc}
 \underline{12} 11 \underline{12}, \\
 \underline{12} 12 \underline{12}, \\
\underline{12} 21  \underline{12}, \\
\underline{12} 22 \underline{12}, 
\end{array}
\right\}
$$
Each element of $E_3^1 $ has 
summands all of which have the same   
first and last factor, and therefore cannot 
appear in $ E_3^2 $. 
Since there are clearly no duplicates inside the two sets, we get then i). 
Applying $ u $ to the two sets produces the same 
two sets, with the underlines removed, so also iii) and iv) follow. 
But then also 
ii) follows: $ u $ picks out the highest summand of the elements, so the matrix relating 
the vectors of $ E_3 $ and those of $ u(E_3 ) $ is lower triangular with 
respect to our order. Note furthermore that $ E_3 $ clearly is a subset of 
$ \varphi_n(B_n ) $ by the description of this before the 
Proposition.

The proof of the induction step $ n-1 \implies n $ goes as follows. 
First of all $i)$ 
is clear from the definitions. 
From $iii)$ at level $ n-1 $ we get that 
$$ 
| E_n^1 | = | u( E_n^1) |\,\,\,   \mbox{  and  }  
| E_n^2 | = | u( E_n^2) | 
$$ 
But then we get
$ iv) \implies iii) \implies ii) $ at level $ n $: the first implication since
$$ 
| E_n | = | E_n^1 \cup E_n^2 | = | E_n^1 | + | E_n^2 | = 
|u( E_n^1) |+ | u(E_n^2) | = $$
$$ |u( E_n^1) \cup  u(E_n^2) | = |u( E_n^1 \cup  E_n^2) | = | u( E_n )
| 
$$ 
The second implication, since once again $ u $ defines a lower triangular matrix 
with respect to the order and since $ E_n \subset \varphi_n ( B_n ) $ by the 
description of $ \varphi_n ( B_n ) $. So let us prove $iv)$. Now 
$$ 
u(E_n^1) \cap u(E_n^2) = 
u(E_n^1) \cap \{ 12 w 12 \, |\, w \in \seq_{2n-4}\{1,2 \} \setminus u( E_{n-2} ) \} $$ 
Consider first 
$$ 
u(E_n^1) \cap \{ 12 w 12 \, |\, w \in \seq_{2n-4}\{1,2 \}  \} 
$$ 
Any element of this intersection is on the form 
$ 1t2 = u(e) $ where $ e \in E_n^1 $. 
But applying $ u $ and the rule (\ref{* rule}) 
to the different elements of $ E_n^1 $ we see that 
$ e= 1x2 $ with $ x \in E_{n-1} $. 
Now $ u(e) $ is also on the form $12w12 $, so $ x= 2y1 $. But such an 
$ x $ must come from $ E_{n-1}^1 $ and thus $ y \in E_{n-2} $. 
All in all: 
$ u(e)= 12y12 $ with $ y \in E_{n-2} $. 
But then our first intersection is empty proving the last part of the claim. 



\subsection{Fullness and more multiplicities}

The standard content of 
the individual $\ourrep$--permutation modules 
$\permmod{\lambda}{n}$ may be determined
similarly to that of $\ourrep$,  
except that we need to recurse all the $\lambda$s together.


By Lemma~\ref{lemma 2} there is a function $v_{\our}^{\lambda}(\mu)$
such that  
$$
v_{\our}^{\lambda}(\mu) = \SCh{ \permmod{\lambda}{n}  : \standard{\mu}{n}  }
$$ 
for any $n$.

By virtue of (\ref{permres1}) and
(\ref{S res}) we have  
\beqa
v_{\our}^{\lambda-2}(\mu) + 2 v_{\our}^{\lambda}(\mu)
+ v_{\our}^{\lambda+2}(\mu) 
   &=&  v_{\our}^{\lambda}(\mu +1) + v_{\our}^{\lambda}(\mu -1) 
\eeqa
As before, explicit inspection of the smallest cases is sufficient to
prime a recursion using this formula to determine all multiplicities. 
We  have for example 
\[
\begin{array}{r|cccccccccc}
&&&&& \lambda = 
\\
&8 & 6 & 4 & 2 & 0 & -2 & ..
\\
\hline
\\
     -4&1 & 7 & 19 & 31 & 37 & 31 & ..
\\
     -3&  & 1 &  5 &  9 & 11 &  9 & 5 & 1
\\
     -2&  &   &  1 &  3 &  3 &  3 & 1
\\
     -1&  &   &    &  1 &  1 &  1 &
\\
\mu = 0&  &   &    &    &  1 &
\\
+1     &  &   &    &    &  1 &
\\
   +2  &  &   &    &  1 &  1 &  1 &
\\
    +3 &  &   &  1 &  3 &  3 &  3 & 1
\\
    +4 &  & 1 &  5 &  9 & 11 &  9 & 5 & 1
\\
    +5 &1 & 7 & 19 & 31 & 37 & 31 & ..
\end{array}
\]
In the format of this table, for every subpart of form 
\eql(croos)
\begin{array}{ccc} &x \\ a&b&c \\ &y \end{array}
\eq
we have $a+2b+c = x+y$. 

 

Note that $\permmod{\lambda}{} \cong \permmod{-\lambda}{}$ 
and that such a $\lambda$ is  necessarily even. 
Accordingly, call $\{ 2n, 2n-2, 2n-4, .. , 0 \}$ the set of
$\Permmod$--weights of $b_n$ --- a sufficient set of labels for
inequivalent permutation modules. 

We do not need a closed formula for all the multiplicities,
but rather  

\prl(prefull1) 
Restrict attention to $\lambda $ an $M$--weight. 
Then
$$  v_{\our}^{\lambda}({\mu}) =   
\left\{ \begin{array}{ll} 
1 & \mbox{$2 \mu = -\lambda$}
\\ 
1 & \mbox{$2 \mu = (\lambda+2)$}
\\ 
0 & \mbox{$0 > 2 \mu > -\lambda$} 
\\ 
0 & \mbox{$0 < 2 \mu  < ( \lambda +2)$} 
\end{array} \right.
$$
\end{pr}
{\em Proof:}
This is the neighbourhood of the domain of zeros (unwritten) in our
table above. The template (\ref{croos}) populates this region as claimed
with the rows at $\mu=0,\pm 1$ as base. \Qed
\newline
(Another proof follows from noting, for example, that 
$\permmod{2n}{n} = \Delta_{n}(-n)$ so 
$( \permmod{\lambda}{n} : \Delta_{n}(\mu) )= 0 $ if $0 > 2 \mu > -\lambda $
and 
$( \permmod{\lambda}{n} : \Delta_{n}(\mu) )= 1 $ if $ 2 \mu = -\lambda $.)


\begin{co}
The module $\ourrep$ is full tilting. 
\end{co}
{\em Proof:} The singleton multiplicities in the expression above give
a bijection between the $\Permmod$--weights and ordinary weights. 
Recall \cite{Donkin98} that each $\tilting{\mu}{}$ contains:
\begin{itemize}
\item
one copy of $\Delta_{\mu}$, and 
\item
no copy of any other 
standard module except having weight higher in the heredity order. 
\end{itemize}
The proposition thus implies that $\permmod{\lambda}{}$ contains 
no $\tilting{\mu}{}$ unless $\mu$  lower than 
(or equal to) the ordinary weight corresponding to
$\lambda$; 
and hence exactly one
copy of the indecomposable tilting module associated to the
corresponding ordinary weight. 
\Qed


\section{On the generic standard content of $\altrep$}


The question of tilting  for $\altrep$ remains open 
(our specific combinatorial constuction in the proof of injectivity of
$\varphi_n$ is particular to $\ourrep$). 
For the reasons outlined in \cite{MartinWoodcock02} it might be useful
to know the standard content of $\altrep$ when it {\em is} tilting. 
Just as for $\ourrep$ 
we have 
\beqa \label{F2}
 F( \altrep (n+2)  )   & \cong &  \altrep (n)
\eeqa
$$ 
\Res(n,n+2,{ \altrep (n+2)}) =4  \altrep (n)
$$
The argument for (\ref{F2}) in the $\altrep$ case is exactly the same as
before. 


Under the assumption that $\altrep$  
has a standard filtration  
(as
in any semisimple specialisation for example), 
it follows 
from Proposition~\ref{blobApply F} and (\ref{F2}) 
that there is a function 
$v': \Z \rightarrow \N$   
such that
$$
\SCh{ \altrep(n) : \standard{\lambda}{}} = v'(\lambda)  
$$ 
(any $n$, $|\lambda|\leq n$, $\lambda -n \equiv 0$ mod.2). 
Let $\altM(i,j):\Z \rightarrow \N$ be
$\altM(i,j)=\delta_{i,j-1}+\delta_{i,j+1}$, 
so  (from \cite{MartinWoodcock2000}) 
\beqa \label{resstan}
\SCh{\Res(n,n+1,\standard{\mu}{n+1}) : \standard{\lambda}{n}} 
  & = &  \altM(\lambda,\mu)
\eeqa
Regarding $\altM,v'$   
as infinite matrices it follows that 
$$ 
\altM v' = 4v' 
$$
which is to say that 
$$ 
v' (\lambda +1) + v'(\lambda -1) = 4 v'(\lambda ) 
$$
Thus $v'$
is determined by recursion from the initial conditions 
$$
v'(0) = 1, \;\;\;\;\; v'(1)=2, \;\;\;\;\; v'(-1)=2 
$$
(which may be determined by inspection of the representations
themselves --- note that it is only these initial conditions which
distinguish this analysis from a corresponding one for $\ourrep$). 
In case $l>0$ we may now obtain $v'(l+1)$   
by 
 $v'(l+1)=4v'(l)-v'(-l+1)=4v'(l)-v'(l-1)$ 
($l<0$ case similar, or note that $v'(-l)=v'(l)$).
We have 
\[
\begin{array}{r|rrrrrrrrrrr}
l     &  4  & 3  & 2  & 1 & 0 & -1 & -2 & -3 \\
v'(l) & 97  & 26 & 7  & 2 & 1 &  2 &  7 & 26
\end{array}
\]


\section{Discussion}

In \cite{MartinWoodcock02} the representations $\ourrep$ and $\altrep$
were introduced, and Martin and Woodcock posed the question of whether
these representations are full tilting. We have now answered this
question in the affirmative for $\ourrep$. 
(They also asked if the representations are faithful for arbitrary $k$
--- a question we answer in the affirmative in \cite{MartinRyom02b}, 
using entirely different techniques.) 
The primary focus of the original paper, however, was 
{\em  generalisations} of the blob algebra. 
In particular it points out the potential
usefulness of corresponding generalisations of $\ourrep$. 
It does not succeed in constructing any. 
The discovery in the present paper that $\ourrep$ is tilting makes it
even more desirable to find such generalisations. 

Since we have constructed a full tilting module for $b_n$ we have,
formally at least, constructed a Ringel dual, 
${\cal B}_n = \mbox{End}_{b_n}(V^{\otimes 2n})$. 
Armed with this mechanism (and the associated combinatorics,
summarized generically by the truncation 
\[
\mat{ccccccccc}
  &   &   &   & 1 &   &   \\
  &   &   & 1 &   & 1 &   \\
  &   & 1 &   & 2 &   & 1 \\
  & 1 &   & 3 &   & 3 &   & 1 \\
1 &   & 4 &   & 6 &   & 4 &   & 1
\tam 
  \mat{c} 41 \\ 11 \\ 3 \\ 1 \\ 1 \\ 3 \\ 11 \\ 41 \\ 153 \tam 
= \mat{c} 1 \\ 4 \\ 16 \\ 64  \\ 256 \tam 
\]
where the $n^{th}$ matrix row gives the dimensions of standards of
$b_{n-1}$; and the column vector gives their multiplicities, and hence the
dimensions of (co)standards of the dual)
we can search for Lie theoretic settings (i.e., a familiar
presentation) for this dual. 

This search will be the subject of a separate paper, 
but it behoves us to assemble the clues which are now ready to hand. 
In particular, let us look briefly at the most interesting case 
in characteristic~0. This means, essentially, $q$ an $l^{th}$ root of
unity and $m$ an integer ($|m|<l$). 
(Although the connection with Lie theory is still, for the present,
`virtual' we know from \cite{MartinWoodcock2000} that Lie theoretic
terminology provides the correct setting for a description of blob representation
theory.) 
Then the alcove structure is as follows. 
The weight space is $\R$ and integral weights $\Z$.  
The affine Weyl group is generated by a reflection at $m$ and another
at $m-l$. No `wall' (reflection point) lies at 0, so call the alcove
containing 0 the 0--alcove. Label the first alcove on the $\pm$--ve
side of the 0--alcove the $\pm 1$--alcove. Label all other alcoves by
the obvious counting scheme. 
The blocks are the affine Weyl orbits, 
and the regular blocks are (up to localisation)
Morita equivalent, so we will pick one arbitrarily and relabel weights
in it simply by their alcove labels. Then the simple submodule
structure of standard $\Delta(\nu)$ ($\nu \geq 0$) is 
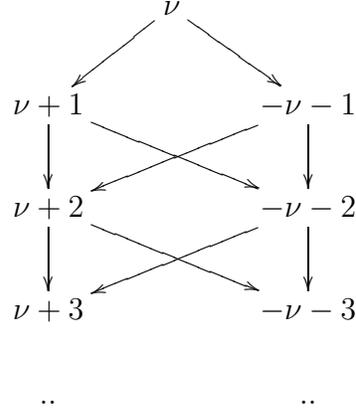
\[ \xymatrix{
& \nu \ar[dl] \ar[dr] \\
\nu +1  \ar[d] \ar[drr] &  & -\nu -1 \ar[dll] \ar[d] \\
\nu +2  \ar[d] \ar[drr] &  & -\nu -2 \ar[dll] \ar[d] \\
\nu +3                  &  & -\nu -3 \\
..                      &  & ..
} \]
(the ladder continues down until truncated by localisation). 
So far all is taken from  \cite{MartinWoodcock2000}.  Now consider
what we may deduce about the indecomposable tilting module labelled by
$\nu$. We have that every simple in the defining standard must be the
socle of a costandard. We need then to take a standard filtered
closure, and assemble the resultant melange into a contravariant
selfdual module. For example: 
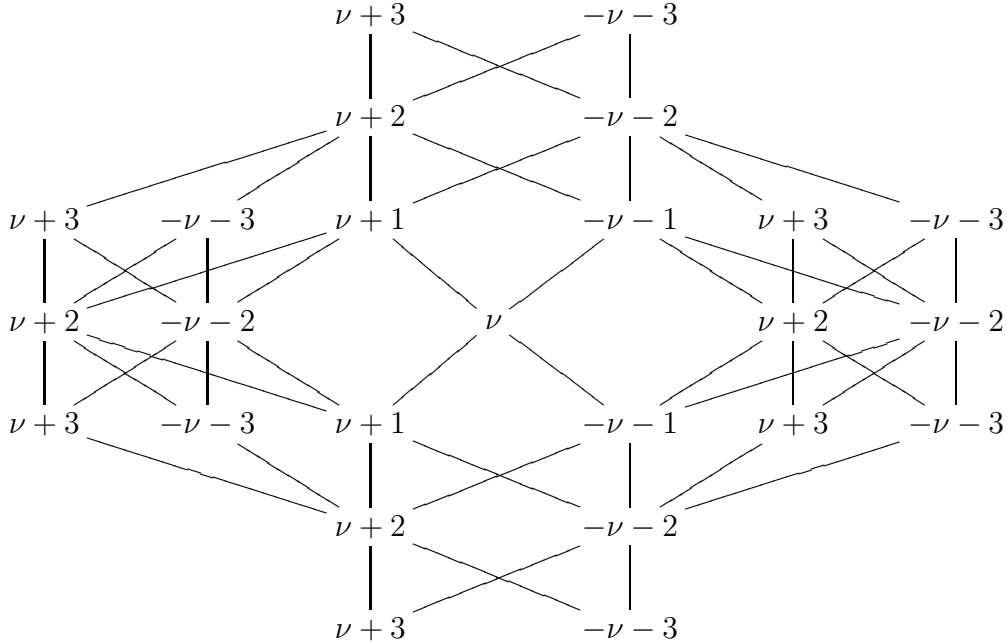
\newcommand{\mybar}{\ar@{-}}
\[ 
 \xymatrix{
&&                  \nu +3  \mybar[d] \mybar[drr] & 
          & -\nu -3 \mybar[dll] \mybar[d] \\
&& \mybar[dl] \mybar[dll] \nu +2  \mybar[d] \mybar[drr] & 
          & -\nu -2 \mybar[dll] \mybar[d] \mybar[dr]\mybar[drr]\\
\nu+3 \mybar[d]\mybar[dr]&-\nu-3 \mybar[d]\mybar[dl]& \nu+1 \mybar[dll]\mybar[dl] \mybar[dr]  & 
          & -\nu -1\mybar[dl]\mybar[drr]\mybar[dr]&\nu+3\mybar[d]\mybar[dr]&-\nu-3\mybar[d]\mybar[dl] \\
\nu+2\mybar[d]\mybar[dr]\mybar[drr]&-\nu-2\mybar[d]\mybar[dl]\mybar[dr]&& 
  \nu \mybar[dl] \mybar[dr] &&\nu+2\mybar[d]\mybar[dr]\mybar[dl]& -\nu-2\mybar[d]\mybar[dl]\mybar[dll] \\
\nu+3\mybar[drr]&-\nu-3\mybar[dr]&\nu +1  \mybar[d] \mybar[drr] & 
          & -\nu -1 \mybar[dll] \mybar[d] &\nu+3\mybar[dl]&-\nu-3 \mybar[dll]\\
&&\nu +2  \mybar[d] \mybar[drr] & 
          & -\nu -2 \mybar[dll] \mybar[d] \\
&&\nu +3                  & 
          & -\nu -3 \\
}
 \]
Here, since layers may contain modules with multiplicity, some of the edges in
the graph indicate no more than layer constraints (although they
provide a useful guide to the eye). These modules are, of course, far
from projective. 

Note that although $T_n(q)$ fails to be \qh\ when $[2]=0$ this failure
is degenerate rather than exceptional, in the sense that if one allows
the notion of a single `formal' standard module of dimension 0 then
the whole formalism is resurrected (the fact that $V^{\otimes n}$
itself is not compromised by passing to $[2]=0$ is a signal of this). 
Similar statements apply in the blob case and, 
as mentioned above, in the paper \cite{MartinRyom02b} we
show that $\ourrep$ is faithful for arbitrary (not just \qh)
specialisations. 
The questions of tilting and faithfulness for $\altrep$ remain open. 


\medskip
\medskip
\noindent
{\bf Acknowledgements.} We would like to thank Anton Cox for several
useful discussions. PPM would like to thank Sheila Brenner, 
Andrew Mathas, and Richard
Green, each for a useful  discussion. We would like to thank EPSRC for
funding under GR/M22536.

\bibliographystyle{amsplain}
\bibliography{new31,main}
\end{document}